\documentclass[11pt]{article}
\usepackage{amssymb,amsthm,amsmath}
\def\qed{\hfill $\Box$}

\newcommand\pf{\smallbreak\noindent \texttt{Proof}. }

\begin{document}

\newtheorem{thm}{Theorem}[section]
\newtheorem{prop}[thm]{Proposition}
\newtheorem{lem}[thm]{Lemma}
\newtheorem{cor}[thm]{Corollary}
\newtheorem{ex}[thm]{Example}
\renewcommand{\thefootnote}{*}

\title{\bf On the automorphism groups of some Leibniz algebras}

\author{\textbf{L.A.~Kurdachenko, A.A.~Pypka}\\
Oles Honchar Dnipro National University, Dnipro, Ukraine\\
{\small e-mail: lkurdachenko@gmail.com, sasha.pypka@gmail.com}\\
\textbf{I.Ya.~Subbotin}\\
National University, Los Angeles, USA\\
{\small e-mail: isubboti@nu.edu}}
\date{}

\maketitle

\begin{abstract}
We study the automorphism groups of finite-dimensional cyclic Leibniz algebras. In this connection, we consider the relationships between groups, modules over associative rings and Leibniz algebras.
\end{abstract}



\thispagestyle{empty}

\section{Introduction}
Let $L$ be an algebra over a field $F$ with the binary operations $+$ and $[,]$. Then $L$ is called a {\em left Leibniz algebra} if it satisfies the {\em left Leibniz identity}
$$[[a,b],c]=[a,[b,c]]-[b,[a,c]]$$
for all $a,b,c\in L$.

Leibniz algebras first appeared in the paper of A.~Blokh~\cite{BA1965}, while the term ``Leibniz algebra'' appeared in the book of J.-L.~Loday~\cite{LJ1992}, and the article of J.-L.~Loday~\cite{LJ1993}. In~\cite{LP1993}, J.-L.~Loday and T.~Pirashvili began to study properties of Leibniz algebras. The theory of Leibniz algebras has developed very intensively in many different directions. Some of the results of this theory were presented in the book~\cite{AOR2020}. Note that Lie algebras are a partial case of Leibniz algebras. Conversely, if $L$ is a Leibniz algebra, in which $[a,a]=0$ for every $a\in L$, then it is a Lie algebra. Thus, Lie algebras can be characterized as anticommutative Leibniz algebras. The question about those properties of Leibniz algebras that are absent in Lie algebras, and accordingly about those types of Leibniz algebras that have essential differences from Lie algebras, naturally arises. A lot has already been done in this direction, including some results of the authors of this article. Many new results can be found in the survey papers~\cite{CPSY2019,KKPS2017,KSeSu2020}. Other results related to this topic can be found in~\cite{AAO2005,AAO2006,BD2011,BD2013,GVKO2012,KOS2019,KSeSu2017,KSuSe2018,KSY2018}.

In the study of Leibniz algebras the information about their endomorphisms and derivations is very useful, as shown, for example, in~\cite{AKOZh2020,LRT2016}. The endomorphisms and derivations of infinite-dimensional cyclic Leibniz algebras were investigated in~\cite{KSY2021}.

Let $L$ be a Leibniz algebra. As usual, a linear transformation $f$ of $L$ is called an {\em endomorphism} of $L$ if
$$f([a,b])=[f(a),f(b)]$$
for all $a,b\in L$. Clearly, a product of two endomorphisms of $L$ is also an endomorphism of $L$, so that the set of all endomorphisms of $L$ is a semigroup by a multiplication. In the same time, the sum of two endomorphisms of $L$ is not necessary to be an endomorphism of $L$, so that we cannot speak about the endomorphism ring of $L$.

Here we will use the term ``{\em semigroup}'' for the set that has an associative binary operation. For a semigroup with an identity element, we will use the term ``{\em monoid}''. Clearly, an identical transformation is an endomorphism of $L$. Therefore, the set $\mathbf{End}_{[,]}(L)$ of all endomorphisms of $L$ is a monoid by a multiplication. As usual, a bijective endomorphism of $L$ is called an {\em automorphism} of $L$.

Let $f$ be an automorphism of $L$. Then the mapping $f^{-1}$ is also an automorphism of $L$. Indeed, let $x,y$ be arbitrary elements of $L$. Then there are elements $u,v\in L$ such that $x=f(u)$, $y=f(v)$ and
$$f^{-1}([x,y])=f^{-1}([f(u),f(v)])=f^{-1}(f[u,v])=[u,v]=[f^{-1}(x),f^{-1}(y)].$$
Thus, the set $\mathbf{Aut}_{[,]}(L)$ of all automorphisms of $L$ is a group by a multiplication.

It should be noted that endomorphisms of Leibniz algebras have hardly been studied. It is also quite unusual that the structure of cyclic Leibniz algebras is described relatively recently in~\cite{CKS2017}. In this paper, we began the study of the structure of the automorphism groups of finite-dimensional cyclic Leibniz algebras.

Let $L$ be a cyclic Leibniz algebra, $L=\langle a\rangle$, and suppose that $L$ has finite dimension over a field $F$. Then there exists a positive integer $n$ such that $L$ has a basis $a_{1},\ldots,a_{n}$ where
$$a_{1}=a, a_{2}=[a_{1},a_{1}], \ldots, a_{n}=[a_{1},a_{n-1}], [a_{1},a_{n}]=\alpha_{2}a_{2}+\ldots+\alpha_{n}a_{n}.$$
Moreover, $[L,L]=\mathbf{Leib}(L)=Fa_{2}+\ldots+Fa_{n}$~\cite{CKS2017}. We fix these designations.

The following types of cyclic Leibniz algebras appear here.

The first case: $[a_{1},a_{n}]=0$. In this case, $L$ is nilpotent, and we will say that $L$ is a {\em cyclic algebra of type} (I). The structure of the automorphism group of a cyclic Leibniz algebra of type (I) is described in Section~1.

Now we need the following concepts. The {\em left} (respectively {\em right}) {\em center} $\zeta^{\mathrm{left}}(L)$ (respectively $\zeta^{\mathrm{right}}(L)$) of a Leibniz algebra $L$ is defined by the following rule:
$$\zeta^{\mathrm{left}}(L)=\{x\in L|\ [x,y]=0\ \mbox{for each }y\in L\}$$
(respectively,
$$\zeta^{\mathrm{right}}(L)=\{x\in L|\ [y,x]=0\ \mbox{for each }y\in L\}).$$

The left center of $L$ is an ideal of $L$, but it is not true for the right center of $L$. Moreover, $\mathbf{Leib}(L)\leqslant\zeta^{\mathrm{left}}(L)$, so that $L/\zeta^{\mathrm{left}}(L)$ is a Lie algebra. The right center of $L$ is a subalgebra of $L$, and in general the left and right centers are different. They can even have different dimensions (see~\cite{KOP2016}).

The {\em center} $\zeta(L)$ of $L$ is defined by the following rule:
$$\zeta(L)=\{x\in L|\ [x,y]=0=[y,x]\ \mbox{for each }y\in L\}.$$
The center is an ideal of $L$.

Consider now the second type of cyclic Leibniz algebras. In this case, $[a_{1},a_{n}]=\alpha_{2}a_{2}+\ldots+\alpha_{n}a_{n}$ and $\alpha_{2}\neq0$. Put $c=\alpha_{2}^{-1}(\alpha_{2}a_{1}+\ldots+\alpha_{n}a_{n-1}-a_{n})$. Then $[c,c]=0$, $Fc=\zeta^{\mathrm{right}}(L)$, $L=[L,L]\oplus Fc$, $[c,b]=[a_{1},b]$ for every $b\in[L,L]$~\cite{CKS2017}. In particular, $a_{3}=[c,a_{2}],\ldots,a_{n}=[c,a_{n-1}],[c,a_{n}]=\alpha_{2}a_{2}+\ldots+\alpha_{n}a_{n}$. In this case, we will say that $L$ is a {\em cyclic algebra of type} (II).

For description of the automorphism group of a cyclic Leibniz algebra of type (II) we consider the relationships between Leibniz algebras and modules over associative ring. We will consider these relationships in Section~2. Using these constructions, in Sections~3 we obtain the description of automorphism group of a cyclic Leibniz algebra of type (II).

The third case: $[a_{1},a_{n}]=\alpha_{2}a_{2}+\ldots+\alpha_{n}a_{n}$ and $\alpha_{2}=0$. Let $t$ be the first index such that $\alpha_{t}\neq0$. In other words, $[a_{1},a_{n}]=\alpha_{t}a_{t}+\ldots+\alpha_{n}a_{n}$. By our condition, $t>2$. Then
$$[a_{1},a_{n}]=\alpha_{t}[a_{1},a_{t-1}]+\ldots+\alpha_{n}[a_{1},a_{n-1}]=[a_{1},\alpha_{t}a_{t-1}+\ldots+\alpha_{n}a_{n-1}],$$
which implies that $\alpha_{t}a_{t-1}+\ldots+\alpha_{n}a_{n-1}-a_{n}\in\mathbf{Ann}_{L}^{\mathrm{right}}(a_{1})$. Since $\alpha_{t}\neq0$, then $\alpha_{t}^{-1}\neq0$ and $d_{t-1}=\alpha_{t}^{-1}(\alpha_{t}a_{t-1}+\ldots+\alpha_{n}a_{n-1}-a_{n})=a_{t-1}+\beta_{t}a_{t}+\ldots+\beta_{n}a_{n}\in\mathbf{Ann}_{L}^{\mathrm{right}}(a_{1})$.
Put
\begin{equation*}
\begin{split}
d_{t-2}&=a_{t-2}+\beta_{t}a_{t-1}+\ldots+\beta_{n}a_{n-1},\\
d_{t-3}&=a_{t-3}+\beta_{t}a_{t-2}+\ldots+\beta_{n}a_{n-2},\\
&\ldots\\
d_{1}&=a_{1}+\beta_{t}a_{2}+\ldots+\beta_{n}a_{n-t+1}.
\end{split}
\end{equation*}
Then
\begin{equation*}
\begin{split}
[d_{1},d_{1}]&=[a_{1},d_{1}]=d_{2},\\
[d_{1},d_{2}]&=[a_{1},d_{2}]=d_{3},\\
&\ldots\\
[d_{1},d_{t-2}]&=[a_{1},d_{t-2}]=d_{t-1},\\
[d_{1},d_{t-1}]&=[a_{1},d_{t-1}]=0.
\end{split}
\end{equation*}
It follows that the subspace $U=Fd_{1}\oplus Fd_{2}\oplus\ldots\oplus Fd_{t-1}$ is a nilpotent subalgebra. Moreover, a subspace $[U,U]=Fd_{2}\oplus\ldots\oplus Fd_{t-1}$ is an ideal of $L$. Put further $d_{t}=a_{t},d_{t+1}=a_{t+1},\ldots,d_{n}=a_{n}$.

The following matrix corresponds to this transaction:
\begin{equation*}
\left(\begin{array}{ccccccccccc}
1 & \beta_{t} & \beta_{t+1} & \ldots & \beta_{n} & 0 & 0 & \ldots & 0 & 0 & 0\\
0 & 1 & \beta_{t} & \ldots & \beta_{n-1} & \beta_{n} & 0 & \ldots & 0 & 0 & 0\\
0 & 0 & 1 & \ldots & \beta_{n-2} & \beta_{n-1} & \beta_{n} & \ldots & 0 & 0 & 0\\
\ldots & \ldots & \ldots & \ldots & \ldots & \ldots & \ldots & \ldots & \ldots & \ldots & \ldots\\
0 & 0 & 0 & \ldots & 1 & \beta_{t} & \beta_{t+1} & \ldots & \beta_{n-1} & \beta_{n} & 0 \\
0 & 0 & 0 & \ldots & 0 & 1 & \beta_{t} & \ldots & \beta_{n-2} & \beta_{n-1} & \beta_{n} \\
0 & 0 & 0 & \ldots & 0 & 0 & 1 & \ldots & 0 & 0 & 0 \\
\ldots & \ldots & \ldots & \ldots & \ldots & \ldots & \ldots & \ldots & \ldots & \ldots \\
0 & 0 & 0 & \ldots & 0 & 0 & 0 & \ldots & 1 & 0 & 0 \\
0 & 0 & 0 & \ldots & 0 & 0 & 0 & \ldots & 0 & 1 & 0 \\
0 & 0 & 0 & \ldots & 0 & 0 & 0 & \ldots & 0 & 0 & 1
\end{array}\right)
\end{equation*}
This matrix is non-singular, which shows that the elements $\{d_{1},\ldots,d_{n}\}$ form a new basis. We note that a subspace $V=Fd_{t}\oplus\ldots\oplus Fd_{n}$ is a subalgebra. Moreover, $V$ is an ideal of $L$, because $[a_{1},d_{t}]=d_{t+1},\ldots,[a_{1},d_{n-1}]=d_{n},[a_{1},d_{n}]=\alpha_{t}d_{t}+\ldots+\alpha_{n}d_{n}$. Moreover, $[a_{1},d_{j}]=[d_{1},d_{j}]$ for all $j\geqslant t$~\cite{CKS2017}. In this case, we will say that $L$ is a {\em cyclic algebra of type}~(III).

Thus, $L=A\oplus Fd_{1}$, $A=V\oplus[U,U]$ is a direct sum of two ideals, $U=[U,U]\oplus Fd_{1}$ is a nilpotent cyclic subalgebra, i.e. is an algebra of type (I), and $V\oplus Fd_{1}$ is a cyclic subalgebra of type (II). The structure of automorphism group of a cyclic Leibniz algebra of type (III) is described in Section~4.

\section{The automorphism group of a cyclic Leibniz algebra of type (I).}
\begin{lem}\label{L1.1}
Let $L$ be a Leibniz algebra over a field $F$, $f$ be an automorphism of $L$. Then $f(\zeta^{\mathrm{left}}(L))=\zeta^{\mathrm{left}}(L)$, $f(\zeta^{\mathrm{right}}(L))=\zeta^{\mathrm{right}}(L)$, $f(\zeta(L))=\zeta(L)$, $f([L,L])=[L,L]$.
\end{lem}
\pf
Let $x$ be an arbitrary element of $L$ and let $z\in\zeta^{\mathrm{left}}(L)$. Since $f$ is an automorphism of $L$, there is an element $y\in L$ such that $x=f(y)$. Then
$$[f(z),x]=[f(z),f(y)]=f([z,y])=f(0)=0.$$
It follows that $f(z)\in\zeta^{\mathrm{left}}(L)$.

On the other hand, there are the elements $u,v\in L$ such that $z=f(u)$, $x=f^{-1}(v)$. Then
$$[u,x]=[f^{-1}(z),f^{-1}(v)]=f^{-1}([z,v])=f^{-1}(0)=0.$$
It follows that $u\in\zeta^{\mathrm{left}}(L)$, so that $z\in f(\zeta^{\mathrm{left}}(L))$ and therefore $f(\zeta^{\mathrm{left}}(L))=\zeta^{\mathrm{left}}(L)$. Using the similar arguments, we obtain that $f(\zeta^{\mathrm{right}}(L))=\zeta^{\mathrm{right}}(L)$ and $f(\zeta(L))=\zeta(L)$.

If $x,y\in L$, then $f([x,y])=[f(x),f(y)]\in[L,L]$. It follows that $f([L,L])\leqslant[L,L]$.

Conversely, let $w\in[L,L]$. Then $w=\alpha_{1}[u_{1},v_{1}]+\ldots+\alpha_{t}[u_{t},v_{t}]$ for some $u_{1},v_{1},\ldots,u_{t},v_{t}\in L$, $\alpha_{1},\ldots,\alpha_{t}\in F$. Since $f$ is an automorphism of $L$, there are the elements $a_{1},b_{1},\ldots,a_{t},b_{t}\in L$ such that $u_{j}=f(a_{j})$, $v_{j}=f(b_{j})$, $1\leqslant j\leqslant t$. Then
\begin{equation*}
\begin{split}
w&=\sum\limits_{1\leqslant j\leqslant t}\alpha_{j}[u_{j},v_{j}]=\sum\limits_{1\leqslant j\leqslant t}\alpha_{j}[f(a_{j}),f(b_{j})]\\
&=\sum\limits_{1\leqslant j\leqslant t}\alpha_{j}f([a_{j},b_{j}])=f\left(\sum\limits_{1\leqslant j\leqslant t}\alpha_{j}[a_{j},b_{j}]\right)\in f([L,L]).
\end{split}
\end{equation*}
It follows that $[L,L]\leqslant f([L,L])$ and hence $[L,L]=f([L,L])$.
\qed

Let $L$ be a Leibniz algebra. Define the {\em lower central series} of $L$
$$L=\gamma_{1}(L)\geqslant\gamma_{2}(L)\geqslant\ldots\gamma_{\alpha}(L)\geqslant\gamma_{\alpha+1}(L)\geqslant\ldots\gamma_{\delta}(L)=\gamma_{\infty}(L)$$
by the following rule: $\gamma_{1}(L)=L$, $\gamma_{2}(L)=[L,L]$ and recursively $\gamma_{\alpha+1}(L)=[L,\gamma_{\alpha}(L)]$ for all ordinals $\alpha$ and $\gamma_{\lambda}(L)=\bigcap_{\mu<\lambda}\gamma_{\mu}(L)$ for the limit ordinals $\lambda$. As usually, we say that a Leibniz algebra $L$ is called {\em nilpotent}, if there exists a positive integer $k$ such that $\gamma_{k}(L)=\langle0\rangle$. More precisely, $L$ is said to be {\em nilpotent of nilpotency class $c$} if $\gamma_{c+1}(L)=\langle0\rangle$, but $\gamma_{c}(L)\neq\langle0\rangle$.

Define the {\em upper central series}
$$\langle0\rangle=\zeta_{0}(L)\leqslant\zeta_{1}(L)\leqslant\ldots\zeta_{\alpha}(L)\leqslant\zeta_{\alpha+1}(L)\leqslant\ldots\zeta_{\tau}(L)=\zeta_{\infty}(L)$$
of a Leibniz algebra $L$ by the following rule: $\zeta_{1}(L)=\zeta(L)$ is the center of $L$, and recursively $\zeta_{\alpha+1}(L)/\zeta_{\alpha}(L)=\zeta(L/\zeta_{\alpha}(L))$ for all ordinals $\alpha$, and $\zeta_{\lambda}(L)=\bigcup_{\mu<\lambda}\zeta_{\mu}(L)$ for the limit ordinals $\lambda$.

\begin{cor}\label{C1.2}
Let $L$ be a Leibniz algebra over a field $F$, $f$ be an automorphism of $L$. Then $f(\zeta_{\alpha}(L))=\zeta_{\alpha}(L)$, $f(\gamma_{\alpha}(L))=\gamma_{\alpha}(L)$ for all ordinals $\alpha$. In particular, $f(\zeta_{\infty}(L))=\zeta_{\infty}(L)$, $f(\gamma_{\infty}(L))=\gamma_{\infty}(L)$.
\end{cor}

Using transfinite induction we can derive it from Lemma~\ref{L1.1}.

\begin{lem}\label{L1.3}
Let $L$ be a Leibniz algebra over a field $F$, $f$ be an endomorphism of $L$. Then $f(\gamma_{\alpha}(L))\leqslant\gamma_{\alpha}(L)$ for all ordinals $\alpha$. In particular, $f(\gamma_{\infty}(L))\leqslant\gamma_{\infty}(L)$.
\end{lem}
\pf
If $x,y\in L$, then $f([x,y])=[f(x),f(y)]\in[L,L]$. It follows that $f([L,L])\leqslant[L,L]$. Suppose that we have already proved that $f(\gamma_{\beta}(L))\leqslant\gamma_{\beta}(L)$ for all ordinals $\beta<\alpha$. If $\alpha$ is a limit ordinal, then $\gamma_{\alpha}(L)=\bigcap\limits_{\beta<\alpha}\gamma_{\beta}(L)$. In this case,
$$f(\gamma_{\alpha}(L))=f\left(\bigcap\limits_{\beta<\alpha}\gamma_{\beta}(L)\right)\leqslant\bigcap\limits_{\beta<\alpha}f(\gamma_{\beta}(L))\leqslant\bigcap\limits_{\beta<\alpha}\gamma_{\beta}(L)=\gamma_{\alpha}(L).$$
Suppose now that $\alpha$ is not a limit ordinal. Then $\alpha-1=\delta$ exists and $\gamma_{\alpha}(L)=[L,\gamma_{\delta}(L)]$. By induction hypothesis, $f(\gamma_{\delta}(L))\leqslant\gamma_{\delta}(L)$. Let $w\in L$, $v\in\gamma_{\delta}(L)$. Then
$$f([w,v])=[f(w),f(v)]\in[L,\gamma_{\delta}(L)]=\gamma_{\alpha}(L).$$
It follows that $f([L,\gamma_{\delta}(L)])=f(\gamma_{\alpha}(L))\leqslant\gamma_{\alpha}(L)$.
\qed

\begin{lem}\label{L1.4}
Let $L$ be a cyclic finite-dimensional Leibniz algebra over a field $F$. Let $S$ be a subset of all endomorphisms $f$ of $L$ such that $f(x)\in[L,L]$ for each $x\in L$. Then $S=\{f|\ f\in\mathbf{End}_{[,]}(L),f^{2}=0\}$ and $S$ is an ideal of $\mathbf{End}_{[,]}(L)$ with zero multiplication $f\circ g=0$ for every $f,g\in S$.
\end{lem}
\pf
Let $f$ be an endomorphism of $L$ such that $f(a_{1})\in[L,L]$. Then we have $f(a_{2})=f([a_{1},a_{1}])=[f(a_{1}),f(a_{1})]$.
An equality $[L,L]=\mathbf{Leib}(L)$ and the fact that $\mathbf{Leib}(L)\leqslant\zeta^{\mathrm{left}}(L)$ shows that $f(a_{2})=0$. Similarly,
\begin{equation*}
\begin{split}
f(a_{3})&=f([a_{1},a_{2}])=[f(a_{1}),f(a_{2})]=0,\\
&\ldots\\
f(a_{n})&=f([a_{1},a_{n-1}])=[f(a_{1}),f(a_{n-1})]=0.\\
\end{split}
\end{equation*}
It follows that $f(y)=0$ for all $y\in[L,L]$. We have $f(a_{1})=\gamma_{2}a_{2}+\ldots+\gamma_{n}a_{n}$, so that
$$f^{2}(a_{1})=f(f(a_{1}))=f(\gamma_{2}a_{2}+\ldots+\gamma_{n}a_{n})=\gamma_{2}f(a_{2})+\ldots+\gamma_{n}f(a_{n})=0,$$
and similarly, $f^{2}(a_{j})=0$ for all $j$, $2\leqslant j\leqslant n$. It follows that $f^{2}(x)=0$ for all $x\in L$. It means that $f^{2}$ is a zero endomorphism.

Conversely, let $f$ be an endomorphism of $L$, $f^{2}=0$ and let $f(a_{1})=\gamma_{1}a_{1}+\gamma_{2}a_{2}+\ldots+\gamma_{n}a_{n}$.
Then
\begin{equation*}
\begin{split}
0&=f^{2}(a_{1})=f(f(a_{1}))=f(\gamma_{1}a_{1}+\gamma_{2}a_{2}+\ldots+\gamma_{n}a_{n})\\
&=\gamma_{1}f(a_{1})+(\gamma_{2}f(a_{2})+\ldots+\gamma_{n}f(a_{n}))\\
&=\gamma_{1}^{2}a_{1}+((\gamma_{1}\gamma_{2}a_{2}+\ldots+\gamma_{1}\gamma_{n}a_{n})+(\gamma_{2}f(a_{2})+\ldots+\gamma_{n}f(a_{n})))\\
&=\gamma_{1}^{2}a_{1}+v
\end{split}
\end{equation*}
where $v\in[L,L]$. Since $Fa_{1}\cap[L,L]=\langle0\rangle$, $f^{2}=0$ implies that $\gamma_{1}^{2}a_{1}=0$ and $v=0$. Thus, $\gamma_{1}^{2}=0$ and $\gamma_{1}=0$. Hence, $S=\{f|\ f\in\mathbf{End}_{[,]}(L),f^{2}=0\}=\{f|\ f\in\mathbf{End}_{[,]}(L),f(x)\in[L,L]\ \mbox{for each }x\in L\}$.

Let $f\in S$ and $g$ be an arbitrary endomorphism of $L$. Then $(f\circ g)(x)=f(g(x))\in[L, L]$. Using Lemma~\ref{L1.3}, we obtain that $(g\circ f)(x)=g(f(x))\in[L,L]$. It follows that $S$ is an ideal of $\mathbf{End}_{[,]}(L)$. Moreover, if $f,g\in S$, then $(f\circ g)(x)=f(g(x))=0$, because $g(x)\in[L,L]$.
\qed

As the first step, we consider the structure of the automorphism group of a nilpotent finite-dimensional cyclic Leibniz algebra $L=\langle a\rangle$. In this case, $[a_{1},a_{n}]=0$.

\begin{lem}\label{L1.5}
Let $L$ be a cyclic Leibniz algebra of type $\mathrm{(I)}$ over a field $F$. Then a linear mapping $f$ is an endomorphism of $L$ if and only if
\begin{equation*}
\begin{split}
f(a_{1})&=\gamma_{1}a_{1}+\gamma_{2}a_{2}+\ldots+\gamma_{n}a_{n},\\
f(a_{2})&=\gamma_{1}^{2}a_{2}+\gamma_{1}\gamma_{2}a_{3}+\ldots+\gamma_{1}\gamma_{n-1}a_{n},\\
f(a_{3})&=\gamma_{1}^{3}a_{3}+\gamma_{1}^{2}\gamma_{2}a_{4}+\ldots+\gamma_{1}^{2}\gamma_{n-2}a_{n},\\
f(a_{4})&=\gamma_{1}^{4}a_{4}+\gamma_{1}^{3}\gamma_{2}a_{5}+\ldots+\gamma_{1}^{3}\gamma_{n-3}a_{n},\\
&\ldots\\
f(a_{n-1})&=\gamma_{1}^{n-1}a_{n-1}+\gamma_{1}^{n-2}\gamma_{2}a_{n},\\
f(a_{n})&=\gamma_{1}^{n}a_{n}.
\end{split}
\end{equation*}
\end{lem}
\pf
Put $L_{1}=Fa_{1}\oplus\ldots\oplus Fa_{n}=L$, $L_{2}=Fa_{2}\oplus\ldots\oplus Fa_{n}$, \ldots, $L_{n-1}=Fa_{n-1}\oplus Fa_{n}$, $L_{n}=Fa_{n}$. Then $\gamma_{1}(L)=L_{1}$, $\gamma_{2}(L)=L_{2}$, \ldots, $\gamma_{n}(L)=L_{n}$ and $\zeta_{1}(L)=L_{n}$, $\zeta_{2}(L)=L_{n-1}$, \ldots, $\zeta_{n}(L)=L_{1}$.

Lemma~\ref{L1.3} shows that $f(L_{j})\leqslant L_{j}$ for all $j$, $2\leqslant j\leqslant n$. Put $f(a_{1})=\sum\limits_{1\leqslant j\leqslant n}\gamma_{j}a_{j}$.
\begin{equation*}
\begin{split}
f(a_{2})&=f([a_{1},a_{1}])=[f(a_{1}),f(a_{1})]=\left[\sum\limits_{1\leqslant j\leqslant n}\gamma_{j}a_{j},\sum\limits_{1\leqslant k\leqslant n}\gamma_{k}a_{k}\right]\\
&=\left[\gamma_{1}a_{1},\sum\limits_{1\leqslant k\leqslant n}\gamma_{k}a_{k}\right]=\gamma_{1}\left(\sum\limits_{1\leqslant k\leqslant n}\gamma_{k}[a_{1},a_{k}]\right)\\
&=\gamma_{1}^{2}a_{2}+\gamma_{1}\gamma_{2}a_{3}+\ldots+\gamma_{1}\gamma_{n-1}a_{n};\\
f(a_{3})&=f([a_{1},a_{2}])=[f(a_{1}),f(a_{2})]\\
&=\left[\sum\limits_{1\leqslant j\leqslant n}\gamma_{j}a_{j},\gamma_{1}^{2}a_{2}+\gamma_{1}\gamma_{2}a_{3}+\ldots+\gamma_{1}\gamma_{n-1}a_{n}\right]\\
&=[\gamma_{1}a_{1},\gamma_{1}^{2}a_{2}+\gamma_{1}\gamma_{2}a_{3}+\ldots+\gamma_{1}\gamma_{n-1}a_{n}]\\
&=\gamma_{1}^{3}a_{3}+\gamma_{1}^{2}\gamma_{2}a_{4}+\ldots+\gamma_{1}^{2}\gamma_{n-2}a_{n};\\
f(a_{4})&=f([a_{1},a_{3}])=[f(a_{1}),f(a_{3})]\\
&=\left[\sum\limits_{1\leqslant j\leqslant n}\gamma_{j}a_{j},\gamma_{1}^{3}a_{3}+\gamma_{1}^{2}\gamma_{2}a_{4}+\ldots+\gamma_{1}^{2}\gamma_{n-2}a_{n}\right]\\
&=[\gamma_{1}a_{1},\gamma_{1}^{3}a_{3}+\gamma_{1}^{2}\gamma_{2}a_{4}+\ldots+\gamma_{1}^{2}\gamma_{n-2}a_{n}]\\
&=\gamma_{1}^{4}a_{4}+\gamma_{1}^{3}\gamma_{2}a_{5}+\ldots+\gamma_{1}^{3}\gamma_{n-3}a_{n};\\
&\ldots\\
f(a_{n-1})&=f([a_{1},a_{n-2}])=[f(a_{1}),f(a_{n-2})]\\
&=\left[\sum\limits_{1\leqslant j\leqslant n}\gamma_{j}a_{j},\gamma_{1}^{n-2}a_{n-2}+\gamma_{1}^{n-3}\gamma_{2}a_{n-1}+\gamma_{1}^{n-3}\gamma_{3}a_{n}\right]\\
&=[\gamma_{1}a_{1},\gamma_{1}^{n-2}a_{n-2}+\gamma_{1}^{n-3}\gamma_{2}a_{n-1}+\gamma_{1}^{n-3}\gamma_{3}a_{n}]\\
&=\gamma_{1}^{n-1}a_{n-1}+\gamma_{1}^{n-2}\gamma_{2}a_{n};\\
f(a_{n})&=f([a_{1},a_{n-1}])=[f(a_{1}),f(a_{n-1})]\\
&=\left[\sum\limits_{1\leqslant j\leqslant n}\gamma_{j}a_{j},\gamma_{1}^{n-1}a_{n-1}+\gamma_{1}^{n-2}\gamma_{2}a_{n}\right]\\
&=[\gamma_{1}a_{1},\gamma_{1}^{n-1}a_{n-1}+\gamma_{1}^{n-2}\gamma_{2}a_{n}]\\
&=\gamma_{1}^{n}a_{n}.
\end{split}
\end{equation*}

Conversely, let $x=\lambda_{1}a_{1}+\lambda_{2}a_{2}+\ldots+\lambda_{n}a_{n}$ and $y=\mu_{1}a_{1}+\mu_{2}a_{2}+\ldots+\mu_{n}a_{n}$ be arbitrary elements of $L$. Suppose that a linear mapping $f$ satisfies the above conditions. Then

\begin{equation*}
\begin{split}
[x,y]&=[\lambda_{1}a_{1}+\lambda_{2}a_{2}+\ldots+\lambda_{n}a_{n},\mu_{1}a_{1}+\mu_{2}a_{2}+\ldots+\mu_{n}a_{n}]\\
&=[\lambda_{1}a_{1},\mu_{1}a_{1}+\mu_{2}a_{2}+\ldots+\mu_{n}a_{n}]\\
&=\lambda_{1}\mu_{1}[a_{1},a_{1}]+\lambda_{1}\mu_{2}[a_{1},a_{2}]+\ldots+\lambda_{1}\mu_{n}[a_{1},a_{n}];\\
f([x,y])&=f(\lambda_{1}\mu_{1}[a_{1},a_{1}]+\lambda_{1}\mu_{2}[a_{1},a_{2}]+\ldots+\lambda_{1}\mu_{n}[a_{1},a_{n}])\\
&=\lambda_{1}\mu_{1}f([a_{1},a_{1}])+\lambda_{1}\mu_{2}f([a_{1},a_{2}])+\ldots+\lambda_{1}\mu_{n}f([a_{1},a_{n}]);\\
[f(x),f(y)]&=[f(\lambda_{1}a_{1}+\lambda_{2}a_{2}+\ldots+\lambda_{n}a_{n}),f(\mu_{1}a_{1}+\mu_{2}a_{2}+\ldots+\mu_{n}a_{n}]\\
&=[\lambda_{1}f(a_{1})+(\lambda_{2}f(a_{2})+\ldots+\lambda_{n}f(a_{n})),\mu_{1}f(a_{1})+\ldots+\mu_{n}f(a_{n})]\\
&=[\lambda_{1}f(a_{1}),\mu_{1}f(a_{1})+\mu_{2}f(a_{2})+\ldots+\mu_{n}f(a_{n})]\\
&+[\lambda_{2}f(a_{2})+\ldots+\lambda_{n}f(a_{n}),\mu_{1}f(a_{1})+\mu_{2}f(a_{2})+\ldots+\mu_{n}f(a_{n})].
\end{split}
\end{equation*}
By our conditions, $f(a_{j})\in L_{2}$ for all $j\geqslant2$, so that $\lambda_{2}f(a_{2})+\ldots+\lambda_{n}f(a_{n})\in L_{2}$ and therefore
$$[\lambda_{2}f(a_{2})+\ldots+\lambda_{n}f(a_{n}),\mu_{1}f(a_{1})+\mu_{2}f(a_{2})+\ldots+\mu_{n}f(a_{n})]=0.$$
Thus,
\begin{equation*}
\begin{split}
[f(x),f(y)])&=[\lambda_{1}f(a_{1}),\mu_{1}f(a_{1})+\mu_{2}f(a_{2})+\ldots+\mu_{n}f(a_{n})]\\
&=\lambda_{1}\mu_{1}[f(a_{1}),f(a_{1})]+\lambda_{1}\mu_{2}[f(a_{1}),f(a_{2})]+\ldots+\lambda_{1}\mu_{n}[f(a_{1}),f(a_{n})].
\end{split}
\end{equation*}
By our conditions,
\begin{equation*}
\begin{split}
[f(a_{1}),f(a_{1})]&=f([a_{1},a_{1}]),\\
[f(a_{1}),f(a_{2})]&=f([a_{1},a_{2}]),\\
&\ldots\\
[f(a_{1}),f(a_{n})]&=f([a_{1},a_{n}]),
\end{split}
\end{equation*}
which implies that $f([x,y])=[f(x),f(y)]$. Hence, $f$ is an endomorphism of a Leibniz algebra $L$.
\qed

\begin{cor}\label{C1.6}
Let $L$ be a cyclic Leibniz algebra of type $\mathrm{(I)}$ over a field $F$. Then $\mathbf{End}_{[,]}(L)$ is isomorphic to a submonoid of $\mathbf{M}_{n}(F)$ consisting of all matrices having the following form
\begin{equation*}
\left(\begin{array}{cccccccc}
\gamma_{1} & 0 & 0 & 0 & \ldots & 0 & 0 & 0\\
\gamma_{2} & \gamma_{1}^{2} & 0 & 0 & \ldots & 0 & 0 & 0\\
\gamma_{3} & \gamma_{1}\gamma_{2} & \gamma_{1}^{3} & 0 & \ldots & 0 & 0 & 0\\
\gamma_{4} & \gamma_{1}\gamma_{3} & \gamma_{1}^{2}\gamma_{2} & \gamma_{1}^{4} & \ldots & 0 & 0 & 0\\
\ldots & \ldots & \ldots & \ldots & \ldots & \ldots & \ldots & \ldots\\
\gamma_{n-2} & \gamma_{1}\gamma_{n-3} & \gamma_{1}^{2}\gamma_{n-4} & \gamma_{1}^{3}\gamma_{n-5} & \ldots & \gamma_{1}^{n-2} & 0 & 0\\
\gamma_{n-1} & \gamma_{1}\gamma_{n-2} & \gamma_{1}^{2}\gamma_{n-3} & \gamma_{1}^{3}\gamma_{n-4} & \ldots & \gamma_{1}^{n-3}\gamma_{2} & \gamma_{1}^{n-1} & 0\\
\gamma_{n} & \gamma_{1}\gamma_{n-1} & \gamma_{1}^{2}\gamma_{n-2} & \gamma_{1}^{3}\gamma_{n-3} & \ldots & \gamma_{1}^{n-3}\gamma_{3} & \gamma_{1}^{n-2}\gamma_{2} & \gamma_{1}^{n}
\end{array}\right).
\end{equation*}
\end{cor}

\begin{cor}\label{C1.7}
Let $L$ be a cyclic Leibniz algebra of type $\mathrm{(I)}$ over a field $F$. Then the automorphism group $\mathbf{Aut}_{[,]}(L)$ is isomorphic to a subgroup $\mathbf{AC}(n)$ of $\mathbf{GL}_{n}(F)$ consisting of all matrices having the following form
\begin{equation*}
\left(\begin{array}{cccccccc}
\gamma_{1} & 0 & 0 & 0 & \ldots & 0 & 0 & 0\\
\gamma_{2} & \gamma_{1}^{2} & 0 & 0 & \ldots & 0 & 0 & 0\\
\gamma_{3} & \gamma_{1}\gamma_{2} & \gamma_{1}^{3} & 0 & \ldots & 0 & 0 & 0\\
\gamma_{4} & \gamma_{1}\gamma_{3} & \gamma_{1}^{2}\gamma_{2} & \gamma_{1}^{4} & \ldots & 0 & 0 & 0\\
\ldots & \ldots & \ldots & \ldots & \ldots & \ldots & \ldots & \ldots\\
\gamma_{n-2} & \gamma_{1}\gamma_{n-3} & \gamma_{1}^{2}\gamma_{n-4} & \gamma_{1}^{3}\gamma_{n-5} & \ldots & \gamma_{1}^{n-2} & 0 & 0\\
\gamma_{n-1} & \gamma_{1}\gamma_{n-2} & \gamma_{1}^{2}\gamma_{n-3} & \gamma_{1}^{3}\gamma_{n-4} & \ldots & \gamma_{1}^{n-3}\gamma_{2} & \gamma_{1}^{n-1} & 0\\
\gamma_{n} & \gamma_{1}\gamma_{n-1} & \gamma_{1}^{2}\gamma_{n-2} & \gamma_{1}^{3}\gamma_{n-3} & \ldots & \gamma_{1}^{n-3}\gamma_{3} & \gamma_{1}^{n-2}\gamma_{2} & \gamma_{1}^{n}
\end{array}\right)
\end{equation*}
where $\gamma_{1}\neq0$.
\end{cor}

\begin{cor}\label{C1.8}
Let $L$ be a cyclic Leibniz algebra of type $\mathrm{(I)}$ over a field $F$. Then a monoid of all endomorphisms of $L$ is a union of an ideal $S=\{f|\ f\in\mathbf{End}_{[,]}(L),f^{2}=0\}$ and an automorphism group $\mathbf{Aut}_{[,]}(L)$. Moreover, $S$ is an ideal with zero multiplication $f\circ g=0$ for every $f,g\in S$.
\end{cor}
\pf
Indeed, consider the arbitrary endomorphism $f$ of $L$ and let
$$f(a_{1})=\gamma_{1}a_{1}+\gamma_{2}a_{2}+\ldots+\gamma_{n}a_{n}.$$
If $\gamma_{1}=0$, then Lemma~\ref{L1.4} shows that $f\in S$. If $\gamma_{1}\neq0$, then Corollaries~\ref{C1.6} and~\ref{C1.7} shows that $f$ is an automorphism of $L$.
\qed

\begin{lem}\label{L1.9}
Let $L$ be a cyclic finite-dimensional Leibniz algebra over a field $F$. Let $G=\mathbf{Aut}_{[,]}(L)$ and $U$ be a subset of all automorphisms $f$ of $L$ such that $f(a_{1})=a_{1}+u$ for some $u\in[L,L]$. Then $U$ is a normal subgroup of $G$ and $G/U$ is isomorphic to a subgroup of the multiplicative group of a field $F$.
\end{lem}
\pf
If $f$ is an arbitrary automorphism of $L$, then $f(a_{1})=\lambda a_{1}+u$ for some $\lambda\in F$, $u\in[L,L]$. Lemma~\ref{L1.4} shows that $\lambda\neq0$. We remark that a coefficient $\lambda$ is uniquely defined. Indeed, suppose that $f(a_{1})=\lambda_{1}a_{1}+u_{1}$ for some $\lambda_{1}\in F$, $u_{1}\in[L,L]$. Then $\lambda a_{1}+u=\lambda_{1}a_{1}+u_{1}$. It follows that
$$(\lambda-\lambda_{1})a_{1}=u_{1}-u.$$
Since $Fa_{1}\cap[L,L]=\langle0\rangle$, we obtain that $\lambda-\lambda_{1}=0$ and $u_{1}-u=0$. Thus, $\lambda=\lambda_{1}$ and $u_{1}=u$.

For the automorphism $f^{-1}$ we have
$$f^{-1}(a_{1})=\sigma a_{1}+w$$
for some $\sigma\in F$, $w\in[L,L]$. Then
\begin{equation*}
\begin{split}
a_{1}&=(f^{-1}\circ f)(a_{1})=f^{-1}(f(a_{1}))=f^{-1}(\lambda a_{1}+u)\\
&=\lambda f^{-1}(a_{1})+f^{-1}(u)=\lambda(\sigma a_{1}+w)+f^{-1}(u)\\
&=\lambda\sigma a_{1}+\lambda w+f^{-1}(u).
\end{split}
\end{equation*}
Using Lemma~\ref{L1.1}, we obtain that $f^{-1}(u)\in[L,L]$, so that $\lambda w+f^{-1}(u)\in[L,L]$. Since $Fa_{1}\cap[L,L]=\langle0\rangle$, we obtain that $\lambda\sigma a_{1}=a_{1}$ and $\lambda w+f^{-1}(u)=0$. Thus, $\sigma=\lambda^{-1},w=-\lambda^{-1}f^{-1}(u)$.

If $f,g\in U$, then $f(a_{1})=a_{1}+u$, $g(a_{1})=a_{1}+v$ for some $u,v\in[L,L]$. Then
$$(f\circ g)(a_{1})=f(g(a_{1}))=f(a_{1}+v)=f(a_{1})+f(v)=a_{1}+u+f(v).$$
Lemma~\ref{L1.1} shows that $f(v)\in[L,L]$. It follows that $f\circ g\in U$.

Let now $f\in U$. As we have seen above, $f^{-1}(a_{1})=a_{1}-f^{-1}(u)$. Using again Lemma~\ref{L1.1}, we obtain that $f^{-1}(u)\in[L,L]$, which means that $f^{-1}\in U$. It follows that $U$ is a subgroup of a group $G$.

Let $h$ be an arbitrary element of $G$ and let again $f\in U$. Then $h(a_{1})=\lambda a_{1}+y$ for some $y\in[L,L]$. By proved above, $h^{-1}(a_{1})=\lambda^{-1}a_{1}-\lambda^{-1}h^{-1}(y)$. Thus,
\begin{equation*}
\begin{split}
(h^{-1}\circ f\circ h)(a_{1})&=h^{-1}(f(h(a_{1})))=h^{-1}(f(\lambda a_{1}+y))\\
&=h^{-1}(\lambda f(a_{1})+f(y))=\lambda h^{-1}(f(a_{1}))+h^{-1}(f(y))\\
&=\lambda h^{-1}(a_{1}+u)+h^{-1}(f(y))=\lambda h^{-1}(a_{1})+\lambda h^{-1}(u)+h^{-1}(f(y))\\
&=\lambda(\lambda^{-1}a_{1}-\lambda^{-1}h^{-1}(y))+\lambda h^{-1}(u)+h^{-1}(f(y))\\
&=a_{1}-h^{-1}(y)+\lambda h^{-1}(u)+h^{-1}(f(y)).
\end{split}
\end{equation*}
Using Lemma~\ref{L1.1}, we obtain that $h^{-1}(y),h^{-1}(u),h^{-1}(f(y))\in[L,L]$. It follows that $h^{-1}\circ f\circ h\in U$, so that $U$ is a normal subgroup of $G$.

Finally, define the mapping $\vartheta:G\rightarrow\mathbf{U}(F)$ by the following rule. Let $f$ be an arbitrary automorphism of $L$, $f(a_{1})=\lambda a_{1}+u$ for some $\lambda\in F$, $u\in[L,L]$. Put $\vartheta(f)=\lambda$. If $h$ is another automorphism of $L$, i.e. $h(a_{1})=\sigma a_{1}+y$ for some $\sigma\in F$, $y\in[L,L]$, then
\begin{equation*}
\begin{split}
(f\circ h)(a_{1})&=f(h(a_{1}))=f(\sigma a_{1}+y)\\
&=\sigma f(a_{1})+f(y)=\sigma(\lambda a_{1}+u)+f(y)\\
&=(\sigma\lambda)a_{1}+\sigma u+f(y)\\
&=(\lambda\sigma)a_{1}+\sigma u+f(y).
\end{split}
\end{equation*}
Lemma~\ref{L1.1} implies that
$$\vartheta(f\circ h)=\lambda\sigma=\vartheta(f)\vartheta(h).$$
Hence, $\vartheta$ is a homomorphism of a group $G$ in $\mathbf{U}(F)$. Clearly, $\mathbf{Ker}(f)=U$.
\qed

\begin{cor}\label{C1.10}
Let $L$ be a cyclic Leibniz algebra of type $\mathrm{(I)}$ over a field $F$, $G=\mathbf{Aut}_{[,]}(L)$. Then $G$ is a semidirect product of a normal subgroup $U$, consisting of all automorphisms $f$ of $L$ such that $f(a_{1})=a_{1}+u$ for some $u\in[L,L]$, and a subgroup $D=\{f|\ f\in\mathbf{Aut}_{[,]}(L),f(a_{1})=\gamma a_{1},0\neq\gamma\in F\}$. Moreover, $D$ is isomorphic to the multiplicative group of a field $F$ and $U$ is isomorphic to a subgroup $\mathbf{UC}(n)$ of $\mathbf{AC}(n)$ consisting of matrices having the following form
$$E+\gamma_{2}\sum\limits_{1\leqslant k\leqslant n-1}E_{k+1,k}+\gamma_{3}\sum\limits_{1\leqslant k\leqslant n-2}E_{k+2,k}+\ldots+\gamma_{n}E_{n,1}.$$
\end{cor}
\pf
Let $f$ be a linear mapping of $L$, having in basis $\{a_{1},a_{2},\ldots,a_{n}\}$ the following matrix
$$\sum\limits_{1\leqslant k\leqslant n}\gamma^{k}E_{k,k}$$
where $0\neq\gamma\in F$. Corollary~\ref{C1.7} shows that $f$ is an automorphism of $L$. Denote by $\mathbf{DmC}(n)$ the subset of $\mathbf{AC}(n)$, consisting of matrices
$$\mathbf{D}(\gamma)=\sum\limits_{1\leqslant k\leqslant n}\gamma^{k}E_{k,k}$$
where $\gamma\neq0$. It is not hard to see that $\mathbf{DmC}(n)$ is a subgroup of $\mathbf{AC}(n)$ and $\mathbf{DmC}(n)\cong D$. Clearly, the mapping $\theta:\mathbf{DmC}(n)\rightarrow\mathbf{U}(F)$ defined by the rule
$$\theta\left(\sum\limits_{1\leqslant k\leqslant n}\gamma^{k}E_{k,k}\right)=\gamma$$
is an isomorphism. It shows that $\mathbf{DmC}(n)$, and hence $D$, is isomorphic to a multiplicative group of a field $F$. In particular, it is abelian.

Consider the set of matrices, having the following form
\begin{equation*}
\left(\begin{array}{ccccccc}
1 & 0 & 0 & 0 & \ldots & 0 & 0\\
\gamma_{2} & 1 & 0 & 0 & \ldots & 0 & 0\\
\gamma_{3} & \gamma_{2} & 1 & 0 & \ldots & 0 & 0\\
\gamma_{4} & \gamma_{3} & \gamma_{2} & 1 & \ldots & 0 & 0\\
\ldots & \ldots & \ldots & \ldots & \ldots & \ldots & \ldots\\
\gamma_{n-1} & \gamma_{n-2} & \gamma_{n-3} & \gamma_{n-4} & \ldots & 1 & 0\\
\gamma_{n} & \gamma_{n-1} & \gamma_{n-2} & \gamma_{n-3} & \ldots & \gamma_{2} & 1
\end{array}\right)
\end{equation*}
Each of these matrices is completely defined by its first column. Therefore, we will denote this matrix by $\mathbf{M}(\gamma_{2},\gamma_{3},\ldots,\gamma_{n})$. Denote the set of all such matrices by $\mathbf{UC}(n)$. Clearly, we can write every matrix from $\mathbf{UC}(n)$ in the following form
$$E+\gamma_{2}\sum\limits_{1\leqslant k\leqslant n-1}E_{k+1,k}+\gamma_{3}\sum\limits_{1\leqslant k\leqslant n-2}E_{k+2,k}+\ldots+\gamma_{n}E_{n,1}.$$
Using Corollary~\ref{C1.7}, we can obtain that the matrix of every automorphism from a subgroup $U$ in basis $\{a_{1},a_{2},\ldots,a_{n}\}$ belong to $\mathbf{UC}(n)$, and conversely. It is not difficult to show that $U\cong\mathbf{UC}(n)$. Also it is not hard to prove that for every matrix $M\in\mathbf{AC}(n)$ we have a decomposition
$$M=\mathbf{M}(\gamma_{2},\gamma_{3},\ldots,\gamma_{n})\mathbf{D}(\gamma_{1}),$$
and we obtain that
$$\mathbf{AC}(n)=\mathbf{UC}(n)\mathbf{DmC}(n).$$
An equality $U\cap D=\langle1\rangle$ is obvious.
\qed

Consider now a polynomial ring $F[X]$. Denote by $R(n)$ the ideal of $F[X]$, generated by the polynomial $X^{n}$. Put $z=X+R(n)$. Then every element of a factor-ring $F[X]/R(n)$ has a form
$$\alpha_{0}+\alpha_{1}z+\alpha_{2}z^{2}+\ldots+\alpha_{n-1}z^{n-1},$$
$\alpha_{0},\alpha_{1},\alpha_{2},\ldots,\alpha_{n-1}\in F$, and this representation is unique. It is possible to show that
$$\mathbf{U}(F[X]/R(n))=\{\alpha_{0}+\alpha_{1}z+\alpha_{2}z^{2}+\ldots+\alpha_{n-1}z^{n-1}|\ \alpha_{0}\neq0\}.$$
Put
$$\mathbf{I}(F[X]/R(n))=\{1+\alpha_{1}z+\alpha_{2}z^{2}+\ldots+\alpha_{n-1}z^{n-1}|\ \alpha_{0},\alpha_{1},\alpha_{2},\ldots,\alpha_{n-1}\in F\}.$$
Then it is not difficult to show that $\mathbf{I}(F[X]/R(n))$ is a subgroup of $\mathbf{U}(F[X]/R(n))$.

\begin{thm}\label{TA}
Let $L$ be a cyclic Leibniz algebra of type $\mathrm{(I)}$ over a field $F$. Then $\mathbf{Aut}_{[,]}(L)$ is a semidirect product of a normal subgroup $U\cong\mathbf{I}(F[X]/R(n))$ and a subgroup $D\cong\mathbf{U}(F)$.
\end{thm}
\pf
Corollary~\ref{C1.10} implies that $G=\mathbf{Aut}_{[,]}(L)$ is a semidirect product of a normal subgroup $U\cong\mathbf{UC}(n)$ and a subgroup $D\cong\mathbf{U}(F)$. Let $\Gamma,\Lambda\in\mathbf{UC}(n)$ where
\begin{equation*}
\begin{split}
\Gamma&=E+\gamma_{1}\sum\limits_{1\leqslant k\leqslant n-1}E_{k+1,k}+\gamma_{2}\sum\limits_{1\leqslant k\leqslant n-2}E_{k+2,k}+\ldots+\gamma_{n-1}E_{n,1},\\
\Lambda&=E+\lambda_{1}\sum\limits_{1\leqslant k\leqslant n-1}E_{k+1,k}+\lambda_{2}\sum\limits_{1\leqslant k\leqslant n-2}E_{k+2,k}+\ldots+\lambda_{n-1}E_{n,1}.
\end{split}
\end{equation*}
Put
$$\Gamma\Lambda=E+\delta_{1}\sum\limits_{1\leqslant k\leqslant n-1}E_{k+1,k}+\delta_{2}\sum\limits_{1\leqslant k\leqslant n-2}E_{k+2,k}+\ldots+\delta_{n-1}E_{n,1}.$$
Since $\Gamma\Lambda\in\mathbf{UC}(n)$, $\Gamma\Lambda$ is completely defined by its first column. We have
\begin{equation*}
\begin{split}
\delta_{1}&=\gamma_{1}+\lambda_{1},\\
\delta_{2}&=\gamma_{2}+\gamma_{1}\lambda_{1}+\lambda_{2},\\
&\ldots\\
\delta_{j}&=\gamma_{j}+\gamma_{j-1}\lambda_{1}+\gamma_{j-2}\lambda_{2}+\ldots+\gamma_{1}\lambda_{j-1}+\lambda_{j},\\
&\ldots\\
\delta_{n-1}&=\gamma_{n-1}+\gamma_{n-2}\lambda_{1}+\gamma_{n-3}\lambda_{2}+\ldots+\gamma_{1}\lambda_{n-2}+\lambda_{n-1}.
\end{split}
\end{equation*}
Taking all this into account, we obtain the following isomorphism. Define a mapping
$$\phi:\mathbf{UC}(n)\rightarrow\mathbf{I}(F[X]/R(n))$$
by the following rule: if $\Gamma\in\mathbf{UC}(n)$, i.e.
$$\Gamma=E+\gamma_{1}\sum\limits_{1\leqslant k\leqslant n-1}E_{k+1,k}+\gamma_{2}\sum\limits_{1\leqslant k\leqslant n-2}E_{k+2,k}+\ldots+\gamma_{n-1}E_{n,1},$$
then put $\phi(\Gamma)=1+\gamma_{1}z+\gamma_{2}z^{2}+\ldots+\gamma_{n-1}z^{n-1}$. By proved above, $\phi(\Gamma\Lambda)=\phi(\Gamma)\phi(\Lambda)$ for every $\Gamma,\Lambda\in\mathbf{UC}(n)$. Clearly, the mapping $\phi$ is bijective, so that $\phi$ is an isomorphism.
\qed

\section{On the relationships between Leibniz algebras and modules over associative rings.}
Let $L$ be a Leibniz algebra over a field $F$ and $A$ be an abelian ideal of $L$. Denote by $\mathbf{End}_{F}(A)$ the set of all linear transformations of $A$. Then $\mathbf{End}_{F}(A)$ is an associative algebra by the operations $+$ and $\circ$. As usual, $\mathbf{End}_{F}(A)$ is a Lie algebra by the operations $+$ and $[,]$ where $[f,g]=f\circ g-g\circ f$ for all $f,g\in\mathbf{End}_{F}(A)$.

Let $u$ be an arbitrary element of $L$. Consider the mapping $\mathbf{l}_{u}:A\rightarrow A$, defined by the rule $\mathbf{l}_{u}(x)=[u,x]$, $x\in A$. For every $u,v\in L$ and $\lambda\in F$ we have
\begin{equation*}
\begin{split}
\mathbf{l}_{u}(x+y)&=[u,x+y]=[u,x]+[u,y]=\mathbf{l}_{u}(x)+\mathbf{l}_{u}(y),\\
\mathbf{l}_{u}(\lambda x)&=[u,\lambda x]=\lambda[u,x]=\lambda\mathbf{l}_{u}(x).
\end{split}
\end{equation*}
Hence, $\mathbf{l}_{u}$ is a linear transformation of $A$. Furthermore, $\beta\mathbf{l}_{u}(x)=\beta[u,x]=[\beta u,x]=\mathbf{l}_{\beta u}(x)$ for every $x\in A$, which implies that $\beta\mathbf{l}_{u}=\mathbf{l}_{\beta u}$. Moreover,
$$(\mathbf{l}_{u}+\mathbf{l}_{v})(x)=\mathbf{l}_{u}(x)+\mathbf{l}_{v}(x)=[u,x]+[v,x]=[u+v,x]=\mathbf{l}_{u+v}(x),$$
which follows that $\mathbf{l}_{u}+\mathbf{l}_{v}=\mathbf{l}_{u+v}$. Consider the mapping $\vartheta:L\rightarrow\mathbf{End}_{F}(A)$, defined by the rule $\vartheta(u)=\mathbf{l}_{u}$, $u\in L$. By above equalities, this mapping is linear. A subspace $\mathbf{Im}(\vartheta)$ is a Lie subalgebra of Lie algebra associated with $\mathbf{End}_{F}(A)$. Denote by $\mathbf{SL}(A)$ the associative subalgebra of $\mathbf{End}_{F}(A)$ generated by $\mathbf{Im}(\vartheta)$. Then the action of $L$ on $A$ can be extended in a natural way to the action of $\mathbf{SL}(A)$. Then $A$ become to a module over associative ring $\mathbf{SL}(A)$. This relationship we will use in a following way.

Let $L$ be a cyclic Leibniz algebra of type~(II). In this case,
$$[a_{1},a_{n}]=\alpha_{2}a_{2}+\ldots+\alpha_{n}a_{n}$$
and $\alpha_{2}\neq0$. Put
$$c=\alpha_{2}^{-1}(\alpha_{2}a_{1}+\ldots+\alpha_{n}a_{n-1}-a_{n}).$$
Then $[c,c]=0$, $Fc=\zeta^{\mathrm{right}}(L)$, $L=[L,L]\oplus Fc$ and $[c,b]=[a_{1},b]$ for every $b\in[L,L]$~\cite{CKS2017}. In particular, $a_{3}=[c,a_{2}]$, \ldots, $a_{n}=[c,a_{n-1}]$, $[c,a_{n}]=\alpha_{2}a_{2}+\ldots+\alpha_{n}a_{n}$.

Put $A=[L,L]$. A linear transformation $\mathbf{l}_{c}:A\rightarrow A$ in basis $\{a_{2},\ldots,a_{n}\}$ has the following matrix
\begin{equation*}
\left(\begin{array}{cccccc}
0 & 0 & 0 & \ldots & 0 & \alpha_{2}\\
1 & 0 & 0 & \ldots & 0 & \alpha_{3}\\
0 & 1 & 0 & \ldots & 0 & \alpha_{4}\\
\ldots & \ldots & \ldots & \ldots & \ldots & \ldots\\
0 & 0 & 0 & \ldots & 0 & \alpha_{n-1}\\
0 & 0 & 0 & \ldots & 1 & \alpha_{n}
\end{array}\right).
\end{equation*}
This matrix is non-degenerate. Hence, $\mathbf{l}_{c}$ is an $F$-automorphism of a linear space $A$. We will consider $A$ as an $F\langle g\rangle$-module where $\langle g\rangle$ is an infinite cyclic group and the action of $g$ on $A$ defined by the following rule: $ga=\mathbf{l}_{c}(a)=[c,a]$ for each element $a\in A$.

Consider now the dual situation. Let $A$ be a vector space over a field $F$ and let $R$ be a subalgebra of an associative algebra $\mathbf{End}_{F}(A)$ of all $F$-endomorphisms of $A$. Then we can consider $R$ as a Lie algebra by the operation
$$[f,g]=f\circ g-g\circ f,$$
$f,g\in R$. Choose in a Lie algebra $R$ some Lie-subalgebra $S$. Put $L=A\oplus S$ and define an operation $[,]$ on $L$ by the following rule:
\begin{equation*}
\begin{split}
[f,g]&=f\circ g-g\circ f\ \mbox{for all }f,g\in S,\\
[a,b]&=0\ \mbox{for all }a,b\in A,\\
[a,f]&=0\ \mbox{for all }a\in A, f\in S,\\
[f,a]&=f(a)\ \mbox{for all }a\in A, f\in S,\\
[a+f,b+g]&=[a,b]+[a,g]+[f,b]+[f,g]=f(b)+[f,g].
\end{split}
\end{equation*}
By such definition, the left center of $L$ includes $A$. Let $x,y,z$ be arbitrary elements of $L$. Then $x=a+f$, $y=b+g$, $z=c+h$ for some $a,b,c\in A$, $f,g,h\in S$. We have
\begin{equation*}
\begin{split}
[x,y]&=[a+f,b+g]=[f,g]+f(b),\\
[y,z]&=[b+g,c+h]=[g,h]+g(c),\\
[x,z]&=[a+f,c+h]=[f,h]+f(c).
\end{split}
\end{equation*}
Then
\begin{equation*}
\begin{split}
[[x,y],z]&=[[f,g]+f(b),c+h]=[f,g](c)+[[f,g],h],\\
[x,[y,z]]&=[a+f,[g,h]+g(c)]=[f,[g,h]]+f(g(c)),\\
[y,[x,z]]&=[b+g,[f,h]+f(c)]=[g,[f,h]]+g(f(c)).
\end{split}
\end{equation*}
Since $S$ is a Lie algebra, $[[f,g],h]=[f,[g,h]]-[g,[f,h]]$, and we obtain that
\begin{equation*}
\begin{split}
[x,[y,z]]-[y,[x,z]]&=[f,[g,h]]+f(g(c))-[g,[f,h]]-g(f(c))\\
&=[f,[g,h]]-[g,[f,h]]+f(g(c))-g(f(c))\\
&=[[f,g],h]+[f,g](c)=[[x,y],z].
\end{split}
\end{equation*}
This shows that $L$ is a Leibniz algebra. If the subalgebra $R$ is commutative, then $S$ as a Lie algebra is abelian. In this case, the right center of $L$ includes $S$.

Let now $A$ be a finite-dimensional vector space over a field $F$ and let $c$ be an $F$-automorphism of $A$. Let $R=F\langle c\rangle$ be an associative subalgebra of $\mathbf{End}_{F}(A)$, generated by the automorphism $c$. This subalgebra is commutative. Therefore, $R$ as a Lie algebra is abelian. Then a subspace $Fc$ is a Lie subalgebra of $R$. Using the above construction, we can construct a Leibniz algebra
$$L=A\oplus Fc.$$
By this way, we come to cyclic Leibniz algebra of type (II).

\section{The automorphism group of a cyclic Leibniz algebra of type (II).}
\begin{lem}\label{L3.1}
Let $L$ be a cyclic Leibniz algebra of type $\mathrm{(II)}$ over a field $F$, $D$ be a centralizer of a subspace $Fc$ in a monoid $\mathbf{End}_{[,]}(L)$. Then $D$ is a submonoid of $\mathbf{End}_{[,]}(L)$. Moreover, $C=D\cap\mathbf{Aut}_{[,]}(L)$ is a normal subgroup of $\mathbf{Aut}_{[,]}(L)$.
\end{lem}
\pf
Indeed, if $f,g$ are two endomorphisms of $L$ such that $f(c)=g(c)=c$, then
$$(f\circ g)(c)=f(g(c))=f(c)=c,$$
so that $f\circ g\in D$. Since the identity mapping of $L$ belong to $D$, $D$ is a submonoid of $\mathbf{End}_{[,]}(L)$.

Let $f$ be an arbitrary element of $C$. Then
$$c=(f^{-1}\circ f)(c)=f^{-1}(f(c))=f^{-1}(c),$$
so that $f^{-1}\in C$.

Let $g$ be an arbitrary automorphism of $L$ and $f$ be an element of $C$. Lemma~\ref{L1.1} shows that $g(c)=\alpha c$ for some $0\neq\alpha\in F$. Then
$$(g^{-1}\circ f\circ g)(c)=g^{-1}(f(g(c)))=g^{-1}(f(\alpha c))=g^{-1}(\alpha c)=\alpha g^{-1}(c)=\alpha\alpha^{-1}c=c.$$
Hence, $C$ is a normal subgroup of $\mathbf{Aut}_{[,]}(L)$.
\qed

If $L$ is a cyclic Leibniz algebra of type (II), then every element of $L$ has a form $a+\alpha c$ where $a\in[L,L]$, $\alpha\in F$, and its presentation in such form is unique.

\begin{lem}\label{L3.2}
Let $L$ be a cyclic Leibniz algebra of type $\mathrm{(II)}$ over a field $F$, $D$ be a centralizer of a subspace $Fc$ in a monoid $\mathbf{End}_{[,]}(L)$. Then $D$ is isomorphic to a multiplicative monoid of a factor-ring $F[X]/\mathbf{a}(X)F[X]$ where $\mathbf{a}(X)=\alpha_{2}+\alpha_{3}X+\ldots+\alpha_{n}X^{n-2}-X^{n-1}$.
\end{lem}
\pf
Put $A=[L,L]$. We can consider $A$ as a module over a polynomial ring $F[X]$, if we define the action of a polynomial $\nu_{0}+\nu_{1}X+\ldots+\nu_{k}X^{k}$ on an arbitrary element $a\in A$ by the following rule:
$$(\nu_{0}+\nu_{1}X+\ldots+\nu_{k}X^{k})a=\nu_{0}a+\nu_{1}\mathbf{l}_{c}(a)+\ldots+\nu_{k}\mathbf{l}_{c}^{k}(a).$$
Since
\begin{equation*}
\begin{split}
a_{3}&=[c,a_{2}]=\mathbf{l}_{c}(a_{2}),\\
a_{4}&=[c,a_{3}]=\mathbf{l}_{c}(a_{3})=\mathbf{l}_{c}(\mathbf{l}_{c}(a_{2}))=\mathbf{l}_{c}^{2}(a_{2}),\\
&\ldots\\
a_{n}&=[c,a_{n-1}]=\mathbf{l}_{c}(a_{n-1})=\mathbf{l}_{c}(\mathbf{l}_{c}^{n-3}(a_{2}))=\mathbf{l}_{c}^{n-2}(a_{2}),
\end{split}
\end{equation*}
and the fact that $\{a_{2},a_{3},\ldots,a_{n}\}$ is a basis of $A$, $A$ becomes a cyclic $F[X]$-module.
$$\mathbf{l}_{c}^{n-1}(a_{2})=\mathbf{l}_{c}(\mathbf{l}_{c}^{n-2}(a_{2}))=\mathbf{l}_{c}(a_{n})=\alpha_{2}+\alpha_{3}\mathbf{l}_{c}(a_{2})+\ldots+\alpha_{n}\mathbf{l}_{c}^{n-2}(a_{2}),$$
so that we can define $\mathbf{l}_{c}^{k}(a_{2})$ (and hence $\mathbf{l}_{c}^{k}(a)$ for arbitrary $a\in A$) for each positive integer $k$.

If $f$ is an endomorphism of $L$, then Lemma~\ref{L1.3} shows that $f(A)\leqslant A$. Define now the mapping $f^{\downarrow}:A\rightarrow A$ by the rule: $f^{\downarrow}(a)=f(a)$ for every $a\in A$. It is not hard to prove that $f^{\downarrow}$ is a linear transformation of a vector space $A$. Suppose now that $f\in D$. Then
$$f^{\downarrow}([c,a])=f([c,a])=[f(c),f(a)]=[c,f(a)]=[c,f^{\downarrow}(a)].$$
On the other hand, $f^{\downarrow}([c,a])=f^{\downarrow}(\mathbf{l}_{c}(a))=f^{\downarrow}(Xa)$ and $[c,f^{\downarrow}(a)]=\mathbf{l}_{c}(f^{\downarrow}(a))=Xf^{\downarrow}(a)$. Thus, we obtain that
$$f^{\downarrow}(Xa)=Xf^{\downarrow}(a).$$
In other words, $f^{\downarrow}$ is an endomorphism of the $F[X]$-module $A$.

Further, $f(a_{2})=\beta_{0}a_{2}+\beta_{1}a_{3}+\beta_{2}a_{4}+\ldots+\beta_{n-2}a_{n}$ for some $\beta_{0},\ldots,\beta_{n-2}\in F$. As we have seen above,
\begin{equation*}
\begin{split}
a_{3}&=[c,a_{2}]=\mathbf{l}_{c}(a_{2})=Xa_{2},\\
a_{4}&=\mathbf{l}_{c}^{2}(a_{2})=X^{2}a_{2},\\
&\ldots\\
a_{n}&=\mathbf{l}_{c}^{n-2}(a_{2})=X^{n-2}a_{2}.
\end{split}
\end{equation*}
Hence, we come to the following presentation
\begin{equation*}
\begin{split}
f(a_{2})&=\beta_{0}a_{2}+\beta_{1}Xa_{2}+\beta_{2}X^{2}a_{2}+\ldots+\beta_{n-2}X^{n-2}a_{2}\\
&=(\beta_{0}+\beta_{1}X+\beta_{2}X^{2}+\ldots+\beta_{n-2}X^{n-2})a_{2}.
\end{split}
\end{equation*}
Put
$$\mathbf{d}_{f}(X)=\beta_{0}+\beta_{1}X+\beta_{2}X^{2}+\ldots+\beta_{n-2}X^{n-2}.$$
Then $f(a_{2})=\mathbf{d}_{f}(X)a_{2}$.

If $a$ is an arbitrary element of $A$, then $a=\sigma_{0}a_{2}+\sigma_{1}a_{3}+\sigma_{2}a_{4}+\ldots+\sigma_{n-2}a_{n}$ for some $\sigma_{0},\ldots,\sigma_{n-2}\in F$. As we have seen above,
\begin{equation*}
\begin{split}
a&=\sigma_{0}a_{2}+\sigma_{1}Xa_{2}+\sigma_{2}X^{2}a_{2}+\ldots+\sigma_{n-2}X^{n-2}a_{2}\\
&=(\sigma_{0}+\sigma_{1}X+\sigma_{2}X^{2}+\ldots+\sigma_{n-2}X^{n-2})a_{2}.
\end{split}
\end{equation*}
Then
\begin{equation*}
\begin{split}
f(a)&=f(\sigma_{0}a_{2}+\sigma_{1}a_{3}+\sigma_{2}a_{4}+\ldots+\sigma_{n-2}a_{n})\\
&=\sigma_{0}f(a_{2})+\sigma_{1}f(a_{3})+\sigma_{2}f(a_{4})+\ldots+\sigma_{n-2}f(a_{n})\\
&=\sigma_{0}f(a_{2})+\sigma_{1}f(Xa_{2})+\sigma_{2}f(X^{2}a_{2})+\ldots+\sigma_{n-2}f(X^{n-2}a_{2})\\
&=\sigma_{0}f(a_{2})+\sigma_{1}Xf(a_{2})+\sigma_{2}X^{2}f(a_{2})+\ldots+\sigma_{n-2}X^{n-2}f(a_{2})\\
&=(\sigma_{0}+\sigma_{1}X+\sigma_{2}X^{2}+\ldots+\sigma_{n-2}X^{n-2})f(a_{2})\\
&=(\sigma_{0}+\sigma_{1}X+\sigma_{2}X^{2}+\ldots+\sigma_{n-2}X^{n-2})\mathbf{d}_{f}(X)a_{2}\\
&=\mathbf{d}_{f}(X)(\sigma_{0}+\sigma_{1}X+\sigma_{2}X^{2}+\ldots+\sigma_{n-2}X^{n-2})a_{2}=\mathbf{d}_{f}(X)a.
\end{split}
\end{equation*}
Thus, we can see that an endomorphism $f$ is defined by the polynomial $\mathbf{d}_{f}(X)$.

Conversely, let $g(X)$ be an arbitrary polynomial.

We define the $F[X]$-endomorphism $\mathbf{s}(g)$ of $A$ by the rule:
$$\mathbf{s}(g)(a)=g(X)a, a\in A.$$
It is not hard to see that if $g(X),r(X)$ are two polynomials, then $\mathbf{s}(g)\circ\mathbf{s}(r)=\mathbf{s}(gr)$.

For every $F[X]$-endomorphism $h$ of $F[X]$-module $A$ define the mapping $h^{\uparrow}:L\rightarrow L$ by the following rule: for arbitrary element $x=a+\gamma c$ of $L$ we put
$$h^{\uparrow}(a+\gamma c)=h(a)+\gamma c.$$
Let $y=b+\sigma c$ be another arbitrary element of $L$. Then
\begin{equation*}
\begin{split}
h^{\uparrow}(x+y)&=h^{\uparrow}(a+\gamma c+b+\sigma c)=h^{\uparrow}(a+b+(\gamma+\sigma)c)\\
&=h(a+b)+(\gamma+\sigma)c=h(a)+h(b)+\gamma c+\sigma c\\
&=h(a)+\gamma c+h(b)+\sigma c=h^{\uparrow}(a+\gamma c)+h^{\uparrow}(b+\sigma c)\\
&=h^{\uparrow}(x)+h^{\uparrow}(y).
\end{split}
\end{equation*}
Let $\mu\in F$, then
\begin{equation*}
\begin{split}
h^{\uparrow}(\mu x)&=h^{\uparrow}(\mu(a+\gamma c))=h^{\uparrow}(\mu a+\mu\gamma c)\\
&=h(\mu a)+\mu\gamma c=\mu h(a)+\mu\gamma c\\
&=\mu(h(a)+\gamma c)=\mu h^{\uparrow}(a+\gamma c)\\
&=\mu h^{\uparrow}(x).
\end{split}
\end{equation*}
Hence, $h^{\uparrow}$ is a linear transformation of a vector space $L$. Furthermore,
\begin{equation*}
\begin{split}
h^{\uparrow}([x,y])&=h^{\uparrow}([a+\gamma c,b+\sigma c])=h^{\uparrow}([\gamma c,b])\\
&=h([\gamma c,b])=h(\gamma[c,b])=\gamma h([c,b])=\gamma h(Xb)\\
&=\gamma Xh(b);\\
[h^{\uparrow}(x),h^{\uparrow}(y)]&=[h^{\uparrow}(a+\gamma c),h^{\uparrow}(b+\sigma c)]\\
&=[h(a)+\gamma c,h(b)+\sigma c]=[\gamma c,h(b)]=\gamma [c,h(b)]\\
&=\gamma Xh(b).
\end{split}
\end{equation*}
Thus, $h^{\uparrow}([x,y])=[h^{\uparrow}(x),h^{\uparrow}(y)]$, so that $h^{\uparrow}$ is an endomorphism of a Leibniz algebra $L$. By definition, $h^{\uparrow}(c)=c$, so that $h^{\uparrow}\in D$. Clearly, if $h_{1},h_{2}$ are different $F[X]$-endomorphism of $A$, then $h_{1}^{\uparrow}\neq h_{2}^{\uparrow}$.

If $f\in D$, then
\begin{equation*}
\begin{split}
f(a+\gamma c)&=f(a)+f(\gamma c)=f(a)+\gamma f(c)\\
&=f(a)+\gamma c=f^{\downarrow}(a)+\gamma c\\
&=(f^{\downarrow})^{\uparrow}(a+\gamma c).
\end{split}
\end{equation*}
In other words, $f=(f^{\downarrow})^{\uparrow}$.

Define now the mapping $\vartheta:F[X]\rightarrow D$ by the rule: $\vartheta(g(X))=\mathbf{s}(g)^{\uparrow}$ for every polynomial $g(X)\in F[X]$. If $f$ is arbitrary element of $D$, then $f^{\downarrow}$ is an $F[X]$-endomorphism of $A$. As we have seen above, this mapping is defined by polynomial $\mathbf{d}_{f^{\downarrow}}(X)$, more precisely, $f^{\downarrow}=\mathbf{s}(\mathbf{d}_{f}(X))$. Then
$$f=(f^{\downarrow})^{\uparrow}=(\mathbf{s}(\mathbf{d}_{f}(X)))^{\uparrow}=\vartheta(\mathbf{d}_{f}(X)),$$
so that a mapping $\vartheta$ is surjective.

Let $g(X),r(X)$ be two polynomials. Recall that $\mathbf{s}(gr)=\mathbf{s}(g)\circ\mathbf{s}(r)$. Let now $h_{1},h_{2}$ be two $F[X]$-endomorphism of $A$. If $x=a+\gamma c$ is an arbitrary element of $L$, then
$$(h_{1}\circ h_{2})^{\uparrow}(a+\gamma c)=(h_{1}\circ h_{2})(a)+\gamma c=h_{1}(h_{2}(a))+\gamma c$$
and
$$(h_{1}^{\uparrow}\circ h_{2}^{\uparrow})(a+\gamma c)=h_{1}^{\uparrow}(h_{2}^{\uparrow}(a+\gamma c))=h_{1}^{\uparrow}(h_{2}(a)+\gamma c)=h_{1}(h_{2}(a))+\gamma c.$$
It proves an equality $(h_{1}\circ h_{2})^{\uparrow}=h_{1}^{\uparrow}\circ h_{2}^{\uparrow}$. Now we have
$$\vartheta(g(X)r(X))=(\mathbf{s}(gr))^{\uparrow}=(\mathbf{s}(g)\circ\mathbf{s}(r))^{\uparrow}=\mathbf{s}(g)^{\uparrow}\circ\mathbf{s}(r)^{\uparrow}=\vartheta(g(X))\vartheta(r(X)).$$
Hence, $\vartheta$ is an epimorphism of a multiplicative monoid $F[X]$ on $D$.

If $g(X)\in\mathbf{Ker}(\vartheta)$, then $f=\vartheta(g(X))$ is an identity automorphism of $L$, i.e. $f(x)=x$ for each $x\in L$. Then $f(a)=f^{\downarrow}(a)=a$ for each $a\in A$. This means that $\mathbf{s}(g)$ is an identity automorphism of $A$, so that $g(X)a=a$ for each $a\in A$. In particular, $g(X)a_{2}=a_{2}$. In other words, $(g(X)-1)a_{2}=0$ and then $g(X)-1\in\mathbf{Ann}_{F[X]}(a_{2})$. In other words, $\mathbf{Ker}(\vartheta)=\mathbf{Ann}_{F[X]}(a_{2})+1$. We note that $\mathbf{Ann}_{F[X]}(a_{2})$ is an ideal of a ring $F[X]$. Now it is not hard to prove that a multiplicative monoid $F[X]/\mathbf{Ker}(\vartheta)$ is isomorphic to the multiplicative monoid of a factor-ring $F[X]/\mathbf{Ann}_{F[X]}(a_{2})$. Finally,
$$X^{n-1}a_{2}=[c,a_{n}]=\alpha_{2}a_{2}+\alpha_{3}a_{3}+\ldots+\alpha_{n}a_{n}=(\alpha_{2}+\alpha_{3}X+\ldots+\alpha_{n}X^{n-2})a_{2},$$
so that $\mathbf{Ann}_{F[X]}(a_{2})=\mathbf{Ann}_{F[X]}(A)=\mathbf{a}(X)F[X]$ where $\mathbf{a}(X)=\alpha_{2}+\alpha_{3}X+\ldots+\alpha_{n}X^{n-2}-X^{n-1}$.
\qed

\begin{thm}\label{TB}
Let $L$ be a cyclic Leibniz algebra of type $\mathrm{(II)}$ over a field $F$. Then $\mathbf{Aut}_{[,]}(L)=G$ includes a normal subgroup $C$, which is isomorphic to
$$\mathbf{U}(F[X]/\mathbf{a}(X)F[X]),$$
where $\mathbf{a}(X)=\alpha_{2}+\alpha_{3}X+\ldots+\alpha_{n}X^{n-2}-X^{n-1}$ such that $G/C$ is isomorphic to a subgroup of a multiplicative group of a field $F$.
\end{thm}
\pf
As in Lemma~\ref{L3.1}, denote by $D$ the centralizer of $Fc$ in $\mathbf{End}_{[,]}(L)$ and let $C=D\cap\mathbf{Aut}_{[,]}(L)$ is a centralizer of $Fc$ in $\mathbf{Aut}_{[,]}(L)$. By Lemma~\ref{L1.1}, $f(Fc)=Fc$ for each $f\in C$, and it follows that $G/C$ is isomorphic to a subgroup of a multiplicative group of a field $F$. An equality $C=D\cap\mathbf{Aut}_{[,]}(L)$ and Lemma~\ref{L3.2} imply that $C$ is isomorphic to $\mathbf{U}(F[X]/\mathbf{a}(X)F[X])$ where
$$\mathbf{a}(X)=\alpha_{2}+\alpha_{3}X+\ldots+\alpha_{n}X^{n-2}-X^{n-1}.$$
\qed

\section{The automorphism group of a cyclic Leibniz algebra of type (III).}
\begin{thm}\label{TC}
Let $L$ be a cyclic Leibniz algebra of type $\mathrm{(III)}$ over a field $F$. Then $\mathbf{Aut}_{[,]}(L)$ is a subdirect product of groups $G_{1}$ and $G_{2}$ where $G_{1}$ is a group described in Theorem~\ref{TA}, $G_{2}$ is a group described in Theorem~\ref{TB}.
\end{thm}
\pf
We have $L=A\oplus Fd_{1}$, $A=V\oplus[U,U]$, $U=Fd_{1}\oplus Fd_{2}\oplus\ldots\oplus Fd_{t-1}$ is a nilpotent cyclic subalgebra, i.e. is an algebra of type (I). Moreover, a subspace $[U,U]=Fd_{2}\oplus\ldots\oplus Fd_{t-1}$ is an ideal of $L$. Furthermore, $V=Fd_{t}\oplus\ldots\oplus Fd_{n}$ is an ideal of $L$, and $[a_{1},d_{j}]=[d_{1},d_{j}]$ for all $j\geqslant t$. In other words, $V\oplus Fd_{1}$ is a cyclic subalgebra of type (II).

Let $G$ be an automorphism group of a Leibniz algebra $L$. Since $L/V\cong U$ is a cyclic nilpotent Leibniz algebra, $G_{1}=G/C_{G}(L/V)$ is a group, which has been described in Theorem~\ref{TA}. Since $L/[U,U]\cong V\oplus Fd_{1}$ is a cyclic Leibniz algebra of second type, $G_{2}=G/C_{G}(L/[U,U])$ is a group, which has been described in Theorem~\ref{TB}.

Let $f\in C_{G}(L/V)\cap C_{G}(L/[U,U])$. Then for each $x\in L$ we have $f(x)=x+a_{1}$ where $a_{1}\in V$ and $f(x)=x+a_{2}$ where $a_{2}\in[U,U]$. It follows that $f(x)-x\in V\cap[U,U]=\langle0\rangle$, so that $f(x)=x$. Thus, $C_{G}(L/V)\cap C_{G}(L/[U,U])=\langle1\rangle$ and Remak's theorem yields the embedding of $G$ into the direct product $G_{1}\times G_{2}$.
\qed

\end{document}